\documentclass[14pt]{amsart}
\usepackage{amsmath,graphicx,amscd}

\theoremstyle{plain}
\newtheorem{theorem}{Theorem}[section]

\theoremstyle{remark}

\theoremstyle{definition}

\def \R {\mathbf{R}}
\def \Z {\mathbf{Z}}
\def \C {\mathbf{C}}

\def \S {\mathbf{S}}

\def\FF{\mathcal{F}}

\numberwithin{equation}{section}
\def\PSL{{\rm PSL}(2,\R)}
\def\SLC{{\rm SL}(2,\C)}

\def\SL{{\rm SL}(2,\R)}

\def\uSL{\widetilde{\SL}}

\def\fg{\pi_1}

\begin{document}
\title[On the integral of  $ \log x\frac{dy}{y}-\log y\frac{dx}{x} $]{ On the integral of  $ \log x\frac{dy}{y}-\log y\frac{dx}{x} $ \\ over the A-polynomial curves }
\author {Vu The Khoi}
\address{ Institute of Mathematics, 18 Hoang Quoc Viet road, 10307, Hanoi, Vietnam}
\email{vtkhoi@math.ac.vn}
\thanks{The author was partially supported by the National Basic Research Program of Vietnam }
\begin{abstract}
{In this note, we study the integral of the 1-form $ \log x\frac{dy}{y}-\log y\frac{dx}{x}$ over certain plane curves defined by A-polynomials of  knots. It is quite surprising that a Chern-Simons type invariant of 3-manifolds, which can be geometrically computed,  may be used to get the exact values of those integrals. The arithmetic nature of 
these integrals is still unknown  at the moment and deserved further investigation.}
\end{abstract}
\maketitle
\section{Introduction} 
\vskip0.1cm
The recent work of D. Boyd and F. Rodriguez-Villegas,  \cite{boyd,boydvil1,boydvil2},  has shown a relationship between the hyperbolic volume of 3-manifolds and the logarithmic Mahler measure of $2$-variable polynomials. The {\it logarithmic Mahler measure} of a non-zero polynomial $P\in \Z[x_1,\cdots, x_n]$ is defined as: $$m(P)= \int_0^1\cdots \int_0^1\log |P(e^{2\pi i\theta_1}, \cdots, e^{2\pi i\theta_n})|d\theta_1 \cdots d\theta_n .$$
 In particular, Boyd and Rodriguez-Villegas show that  if $P$ is the A-polynomial associated to a knot  then, in certain cases, $\pi m(P)$ equals the hyperbolic volume of the knot complement. Consequently, by using a classical formula of Humbert, they show that for arithmetic hyperbolic manifolds, up to multiplication by a  known constant, $m(P)$ equals the value of the Dedekind zeta function of a certain imaginary quadratic extension.       

The above results show the rich arithmetic nature of the hyperbolic volume. It is a well-known  philosophy due to W. Thurston that the  volume and Chern-Simons invariant of a hyperbolic 3-manifold correspond to the real and imagine parts of a holomorphic function. Therefore, it is natural to expect that the Chern-Simons  invariant should give interesting arithmetic consequences. 

In this short note, we show that a Chern-Simons type invariant gives us the exact value of the integral of $ \log x\frac{dy}{y}-\log y\frac{dx}{x} $ over a certain path in the A-polynomial curve. In the light of Thurston's philosophy, this integral should be related to a certain unknown imaginary counterpart of the Mahler measure.   

 The paper  is organized as follows, in the next section we will briefly review about  the  A-polynomials and then give an exposition on  the Godbillon-Vey invariant, a kind of Chern-Simons invariant,  and a Schl\"afli-type formula for the Godbillon-Vey invariant. The last section contains our main result, where we get  the exact values of the integral of the 1-form $ \log x\frac{dy}{y}-\log y\frac{dx}{x} $ over the real part of the A-polynomials curve of certain hyperbolic knots.
 
 The author would like to thank the referee for pointing out several inaccuracies in the previous version of this paper.
\section{ The  A-polynomials and the Godbillon-Vey invariant }
\vskip0.1cm
Let $X$ be a manifold whose boundary $\partial X$ is a torus $T.$ For example, if $K\subset M$ be a knot in a 3-manifold then $X$ may be taken as $M - \eta(K),$ where $\eta(K)$ is a tubular neighborhood of $K.$ 
Let $\fg(X)$ be the fundamental group of $X.$ The fundamental group of $\partial X,$ $\fg(T) = \Z\oplus \Z$ is called the \textit{peripheral subgroup}. Two simple, closed curves $\mu, \lambda$ on $T$ which intersect in one point  generate $\fg(T)$ and are called the \textit{meridian} and \textit{longitude} of $X.$ The pair $(\fg(T)\subset \fg(X))$ is a very powerful invariant of $X,$ however it is difficult to work with the group alone. A classical way to study  groups is to look at their representations into a linear space. 

For a 3-manifold group, the space of representations into $\SLC$ was first systematically studied  by Culler and Shalen  \cite{culler} and has found spectacular applications in geometry and topology of 3-manifolds.        

The $\SLC$ character variety of  a knot complement $X$ is defined by:
 $$\chi(X)=\text{Hom}(\fg(X), \SLC)/ \sim .$$
Here, we use the algebro-geometric quotient in which $\rho,\rho' \in \text{Hom}(\fg(X), \SLC)$ and $\rho\sim \rho'$ if and only if $tr(\rho(g)) = tr(\rho'(g))$ for all $g\in \fg(X).$

In general, the character variety is a affine algebraic set of high dimension and therefore it is hard to work with directly.
In \cite{cooper}, a more manageable algebraic set is introduced by projecting $\chi(X)$ into $\C^2$ using the eigenvalues of the meridian and longitude. This algebraic set is defined by a integral polynomial in two variables.

More precisely, let $i^*: \chi(X)\rightarrow \chi(\partial X)$ be the restriction map induced by the inclusion $\fg(\partial X)\rightarrow \fg(X)$ and $t:\C^*\times \C^*\rightarrow \chi(\partial X)$ defined by associating to a pair $(x,y)$ a character $\rho$ such that $$\rho(\mu) =  \left(\begin{array}{cc}x&0 \\0&1/x\end{array}\right), \ \text{and} \  \rho(\lambda) =  \left(\begin{array}{cc}y& 0\\0&1/y\end{array}\right).$$ Let $Z$ be the union of all the connected component  $C$ of $\chi(X)$ such that   $i^*(C)$ is $1$-dimensional then the closure of  $t^{-1}(i^*(Z)) \subset \C\times \C$ is a plane curve.  
The defining equation of this affine curve, up to multiplication by a non-zero constant, is an integral polynomial $A_X(x,y)$ called the {\it A-polynomial } of $X.$ Thus the A-polynomial parameterizes the restriction of the character variety to the peripheral subgroup.

In general, the A-polynomial has a factor of $(y-1)$ which corresponds to the abelian representations. Without ambiguity, we will take the A-polynomials to be the normalized one by dividing out this factor.  

In the following, we list some basic properties of the A-polynomials. 

- The A-polynomial detects the unknot (\cite{boyer,dunfield}) 

-  Let  $X$ be  the complement of a knot in a homology sphere, then $A_X(x,y)$ involves only even power of $x.$

-  The A-polynomial is reciprocal, i.e., $A_X(x,y) = A_X(1/x, 1/y)$ up to some power of $x$ and $y.$

The reader can consult \cite{boyer, cooper,cooperlong, culler, dunfield} for more details about character variety, A-polynomial and their applications  in topology.
 
We recall the notion of the Godbillon-Vey class of a codimension 1 foliation. Let $\FF$ be a codimension 1 foliation on a manifold $M$ which is defined by a 1-form $\tau.$ The Frobenius integrability condition asserts that there is a 1-form $\theta$ such that $d\tau=\theta \wedge \tau.$ Godbillon and Vey \cite{god} observe that the form $\theta\wedge d\theta$ is closed and that its cohomology class $[\theta\wedge d\theta] \in H^3(M;\R)$ depends only on the foliation $\FF$ and is a cobordism invariant of $\FF.$  

Given an oriented closed 3-manifold $M$ and a representation
$\rho:\fg(M)\longrightarrow\PSL$, we can associate to
$\rho$ a flat connection on the principal $\PSL$ bundle over
$M$, $P_{\rho}={\widetilde M}\times_{\rho}\PSL.$
As
$\PSL$ acts on $S^1,$ the flat connection on $P_{\rho}$ induces a
flat connection on the $S^1$-bundle over $M$,
$E_{\rho}={\widetilde M}\times_{\rho}S^1.$
Let $\FF$ be the codimension-one foliation determined by the horizontal
distribution of the flat connection on
$E_{\rho}$ and  $GV(\FF)$ be its Godbillon-Vey class. Under the assumption that the Euler
class of $E_{\rho}$ is torsion, Brooks-Goldman show that we can obtain from  $GV(\FF)$, in a natural way, a 3-form on $M.$ 
The {\it Godbillon-Vey invariant} of $\rho,$ denoted by $GV(\rho)$, is defined to be the integral of this 3-form over $M.$

The Seifert volume of a 3-manifold $M$, denoted by $\{M\}$, was first defined by R. Brooks and W. Goldman in 1984 as follows: 
$$\{M\}=\max \{\ |GV(\rho)|\quad |\quad \rho : \fg (M) \longrightarrow \PSL, e(\rho)\  \text{is torsion}\ \}.$$      
For a  Seifert fibered manifold, the Seifert volume can be computed from its Seifert data as follows.

Let $M\longrightarrow F$ be a closed Seifert fibered manifold over a surface $F.$ If  $M\longrightarrow F$ has $r$
singular fibers then the Seifert data of $M$ is given by   
$$(g ; (p_1, q_1), (p_2, q_2), \cdots, (p_r, q_r))$$
Where $g=$ genus $(F)$ and $(p_i, q_i)$ are integers encoding the singular type of the $i^{th}$ singular  fiber.
Define the {\it Euler number} and the {\it Euler characteristic} of $M\longrightarrow F$ respectively by
 $$e(M\longrightarrow F)=-\sum_1^r \frac{q_i}{p_i},\qquad \chi(M\longrightarrow
F)=2-2g-\sum_1^r\frac{p_1-1}{p_i}.$$
It is well-known that when $e(M\longrightarrow F)=0$ or $\chi(M\longrightarrow F)\ge 0$ then $M$ admits one
of the five Seifert fibered geometries other  than the $\uSL$ geometry. In these cases, \cite{bg2} shows that $\{M\}=0.$

In the case that $\chi(M\longrightarrow F)< 0$ and $e(M\longrightarrow F)\ne 0 $ we know that $M$ admits the
$\uSL$ geometry.  It follows from \cite{bg2} that   $\{M\}$ equals  the volume of $M$ in this geometry which is
given by 
\begin{equation}\label{vsf}
\frac{4\pi^2\chi(M\longrightarrow F)^2}{|e(M\longrightarrow F)|}.  
\end{equation}

 For more details about the Seifert volume see \cite{bg1,bg2}. 
The main results of   \cite{vu} is the development of a cut-and-paste method for computing the Godbillon-Vey invariant and Seifert volume. The cut-and-paste method is originated from gauge theory (see \cite{kk1,kk2}). Its main steps consist of the followings:

-	Interpret the Godbillon-Vey invariant of a representation as a kind of Chern-Simons invariant associated to the universal covering group $\uSL.$

-	Define the Godbillon-Vey invariant on a manifold with boundary $X$ by gauge-fixing a normal form of the flat connection near $\partial X.$

-	On $X,$ we prove a formula which expresses the difference between the Godbillon-Vey invariants of two representations in a family in terms of the boundary holonomies.

- To apply the formula in the previous step to compute the Godbillon-Vey invariant of
surgery manifolds,  for each representation $\rho$ we have to find a path which 
connects $\rho |_{X}$ to a a representation whose Godbillon-Vey invariant is
already known. 

In this paper, we will need the following result.
\begin{theorem}(\cite{vu}, Theorem  5.1(c))  Suppose that $A_t$ is a normal form purely hyperbolic path of flat connections on a manifold with toral boundary $X.$ Let $\rho_t : \fg (X) \longrightarrow \uSL,$ be the corresponding path of holonomies. Denote by $\mu$ and $\lambda$ the generators of $\fg(\partial X)\cong \Z\oplus \Z.$ If $\rho_t(\mu)=(\tanh a(t), k\pi)$ and    $\rho_t(\lambda)=(\tanh b(t), l\pi)$ then :
$$GV(\rho_1) - GV(\rho_0)=4\int_0^1 (\dot ab-a\dot b)dt.$$
\end{theorem} 
Some explanations about the terminology are needed here.  The normal form of a flat connection is a nice form of the connection near the torus boundary obtained by using gauge transformation. A path of flat connections is called {\it purely hyperbolic} if its boundary holonomies $\rho_t(\mu)$ and $\rho_t(\lambda)$ are hyperbolic elements of $\uSL$ for all $t.$

Here we work with the group $\uSL$ instead of $\PSL$ since $\uSL$ is simply connected and therefore we can trivialize all the principal bundle. The group $\uSL$  can be
described as  $\uSL= \left\{(\gamma,\omega)\vert\quad |\gamma|<1, -\infty<w<\infty\right\}$ with
the group operation defined by:

$(\gamma, \omega)(\gamma', \omega')=(\gamma'', \omega'')$ where
\begin{eqnarray}\label{rule1}
\gamma''&=&(\gamma+\gamma'e^{-2i\omega})(1+\bar
\gamma\gamma'e^{-2i\omega})^{-1}\\\label{rule2}
\omega''&=&\omega+\omega'+\frac{1}{2i}\log\{(1+\bar
\gamma\gamma'e^{-2i\omega})(1+
\gamma\bar\gamma'e^{2i\omega})^{-1}  \}. 
\end{eqnarray}          
Here $\log z$ is defined by its principal value. We then call an element
of $\uSL$ {\it elliptic, parabolic} or {\it hyperbolic} if it
covers an element of the corresponding type in $\SL .$ 

 The  $\uSL$ character variety of a knot complement $X$ is defined by  $$\chi_{\uSL}(X):=\text{Hom}(\fg(X), \uSL)/\uSL.$$
 As the knot complement $X$ satisfies $H^2(X, \Z)=0,$ any $\PSL$  representation of $\fg(X)$ can be lifted to an
$\uSL$ representation.  In fact,  it was shown in \cite{vu} section 6,  that the $\uSL$ character variety of a knot complement is basically a periodic family of $\PSL$ character variety which is  part of the real component of the A-polynomial curve.     The reader should consult \cite{vu} for more details about the Godbillon-Vey invariant and method to compute it..
\section{Values of some logarithmic integrals}

In this section we will consider some examples of 2-bridge knots whose exceptional Dehn surgeries yield Seifert fibered manifolds. We will find the Seifert volume of the surgery manifolds in two ways.
On one hand the Seifert volume is computed by the integral of $ \log x\frac{dy}{y}-\log y\frac{dx}{x} $ over certain paths of the A-polynomial curve by using Theorem 2.1 above. On the other hand for a Seifert fibered manifold $M$ which admits the $\uSL$-geometry  the Seifert volume equals the geometric volume given by formula (2.1). Thus we  get the value of the integrals. The reader can consult \cite{vu} section 7  for a detail study of these examples. For the results on Dehn surgeries see \cite{ble, wu}.  

Denote by $X_{p/q}$ the result of $p/q$-surgery on a knot $K\subset \S^3.$ To compute the Godbillon-Vey invariant on $X_{p/q},$ we write $X_{p/q}=X\cup S,$ where $X$ is the knot complement and $S$ is the solid torus. As \cite{vu} Corollary 6.2 tells us that on the solid torus the Godbillon-Vey invariant vanishes, we only need to do the computation on the knot complement $X$. 

{\bf Figure eight knot}
The A-polynomial curve of the figure eight knot is 
$$A(x,y)= y+y^{-1}-(x^4+x^{-4})+x^2+x^{-2}+2.$$ 
Note that we write here the A-polynomial using negative power of $x$ and $y$ to show the symmetric form of it.  

Computation in \cite{vu} shows that  the image  $i^*(\chi_{\uSL}(X))$ in $\chi_{\uSL}(\partial X)$ of the character variety  can be parameterized as follows :
$$ \rho(\mu)=(\tanh a(s), k\pi), k\in \Z \ \text{and}\  \rho(\lambda)=(\tanh b(s), 0),$$
where $a(s)$ and  $b(s)$ are  given by  
$$\cosh(2a)=\frac{s^2+3s+3}{2(s+1)}\ \text{and} \ \cosh (b)=\frac{s^4+5s^3+7s^2+4s+2}{2(s+1)^2},\qquad
s\ge 0.$$ 
Note that the covering map $\uSL \rightarrow \SL$ sends   $ \rho(\mu)=(\tanh a(s), k\pi)$  and  $\rho(\lambda)=(\tanh b(s), 0)$ to
$\left(\begin{array}{cc}
            e^{a(s)} &0\\
            0&e^{-a(s)}
            \end{array} \right)$
and $\left(\begin{array}{cc}
            e^{b(s)} &0\\
            0&e^{-b(s)}
            \end{array} \right)$ respectively.

Under this covering map, the curve $i^*(\chi_{\uSL}(X))$ above corresponds to a curve $C$ lying in the real part of the A-polynomial curve  $A(x,y)=0$  (see figure 1 below for the plot of the curve $C$). Moreover the form $(\dot ab-a\dot b)dt,$ in Theorem 2.1, descends to the form $ \log y\frac{dx}{x} - \log x\frac{dy}{y}$ on the A-polynomial curve. 

Recall that a
character $\rho \in \chi(X)$ extends to $\chi(X_{p/q})$ if and only if $\rho(p\mu + q\lambda)=1.$
Therefore, the points $(x,y)$ on the A-polynomial curve such that $x^py^q=1$ correspond to the characters in $\chi(X_{p/q}).$           
  
It is well-known that $0$-surgery on the figure eight knot gives a torus bundle. As the torus bundle admits a self map
of degree bigger than 1, it follows from \cite{bg2} that $\{X_0\}=0.$ To find the points on the A-polynomial curve corresponding to characters in  $\chi(X_0),$ we solve the equations $A(x,y)=0, y=1.$ The result is that there is a unique character $\rho_0$  represented by the point $P_0=(\frac{1+\sqrt 5}{2}, 1)$ on the curve $C$ in Figure 1. 
\begin{center}
\includegraphics[scale=0.3, angle=-90]{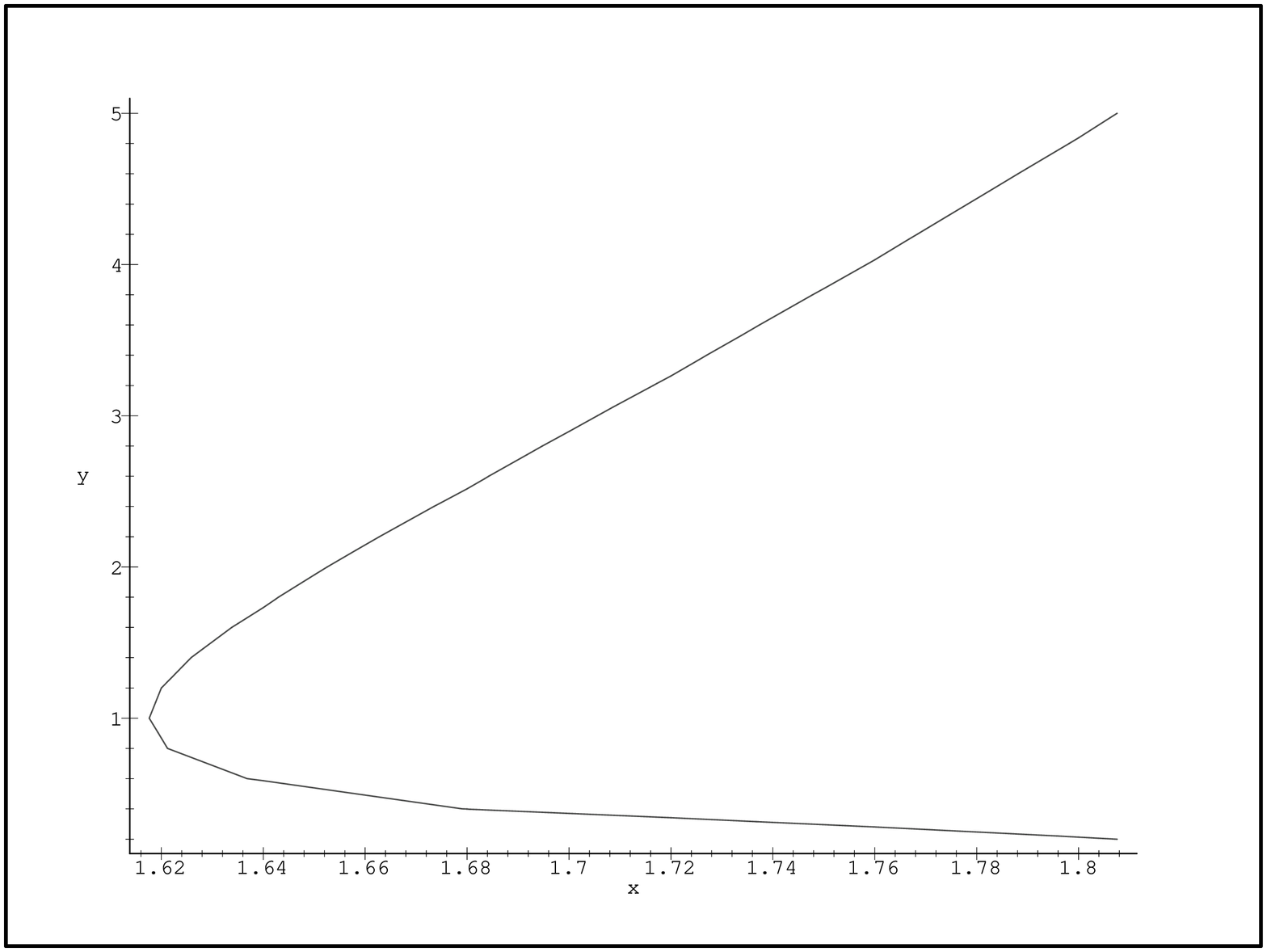}

Fig. 1.  The curve $C$ lying in the A-polynomial curve of the figure eight knot   
\end{center}
{\it (-1)- surgery.} It is known that (-1)-surgery on  the figure eight knot gives the homology sphere $\Sigma(2,3,7).$ This homology sphere carries the $\uSL$ geometry  whose volume is computed to be  $\frac{2\pi^2}{21}$ by using formula (2.1).  The intersection point between the curve $C$ and the line $y/x=1$ corresponds to a representation $\rho_1$ of $\chi(\Sigma(2,3,7)).$ The coordinate of this point is $P_1=(\alpha, \alpha),$ where  $\alpha\approx 1.635573130.$             

Now using Theorem 2.1 we get:
$$GV(\rho_1|_{X}) - GV(\rho_0|{X})=4\int_{\frac{1+\sqrt 5}{2}}^\alpha \log y\frac{dx}{x}-\log x\frac{dy}{y}.$$
where the integral is taken over the curve $C$ from $P_0$ to $P_1$. We have noticed earlier that the contribution from the solid torus is zero, hence  $GV(\rho_1)=GV(\rho_1|_X).$ 
On the other hand, since   $\{X_0\}=0,$ we find that $GV(\rho_0|_X)=0.$ Moreover, since $\rho_1$ is the unique character in $\chi(X_{-1}),$ we get that $\{X_{-1}\}=|GV(\rho_1)|=2\frac{\pi^2}{21}.$ Consequently, we obtained the identity:
$$\int_{\frac{1+\sqrt 5}{2}}^\alpha \log x\frac{dy}{y}-\log y\frac{dx}{x}=\frac{\pi^2}{42}.$$
{\it (-2)-surgery.} For (-2)-surgery , the resulting manifold is the Seifert fibered space over $\S^2$ with three exceptional fibers of indices 2, 4 and 5. We find that $X_{-2}$ admits the $\uSL$ geometry with the volume equals $\frac{\pi^2}{5}.$ The unique intersection point between the curve $C$ and the curve  $x^{-2}y=1$ corresponds to a representation of the surgery manifold. We compute this  intersection point to be $P_2=(\beta, \beta^2)$ where $\beta\approx 1.700015776.$  Arguing similarly as in the (-1)-surgery case, we get :              
  $$\int_{\frac{1+\sqrt 5}{2}}^\beta \log x\frac{dy}{y}-\log y\frac{dx}{x}=\frac{\pi^2}{20}.$$
Here the integral is taken over the curve $C$ from $P_0$ to $P_2.$

{\bf The $5_2$  knot.} This knot is the (7,3) 2-bridge knot  and is indexed by  $5_2$ in the knot table.  Its A-polynomial is given by:
$$A(x,y) = 1+ y(-1+2x^2+2x^4-x^8+x^{10})+y^2(x^4-x^6+2x^{10}+2x^{12}-x^{14})+y^3x^{14}.$$   
Computation from \cite{vu} shows that in this case, the image $i^*(\chi_{\uSL}(X))$ in $\chi_{\uSL}(\partial X)$ of the character variety  consists of three parts. We only interest in the following curve $C$ in $i^*(\chi_{\uSL}(X))$ since it contains characters which extend to characters on the surgery manifolds:
$$\rho(\mu)=(\tanh a(s), k\pi),\ a<0, \ k\in \Z  \  \text{and}\ \rho(\lambda)=(\tanh b(s),0),\
b>0.$$ Here $a(s)$ and  $b(s)$ are given by
\begin{eqnarray*} 
\cosh (2a)&=& \frac{4s^2+6s+4+2\sqrt{s^2+4}}{8s}\\
\cosh (b)&=& 1+\frac{s^5}{4} +s^4+\frac{3s^3}{2}+\frac{3s^2}{2}-(\frac{s^4}{4}s^3+s^2)\sqrt{s^2+4},
\ s>0.
\end{eqnarray*}
 We plot the curve $C$  as part of the real part of the A-polynomial curve
 in Figure. 2 below.
\begin{center}
\includegraphics[scale=0.3, angle=-90]{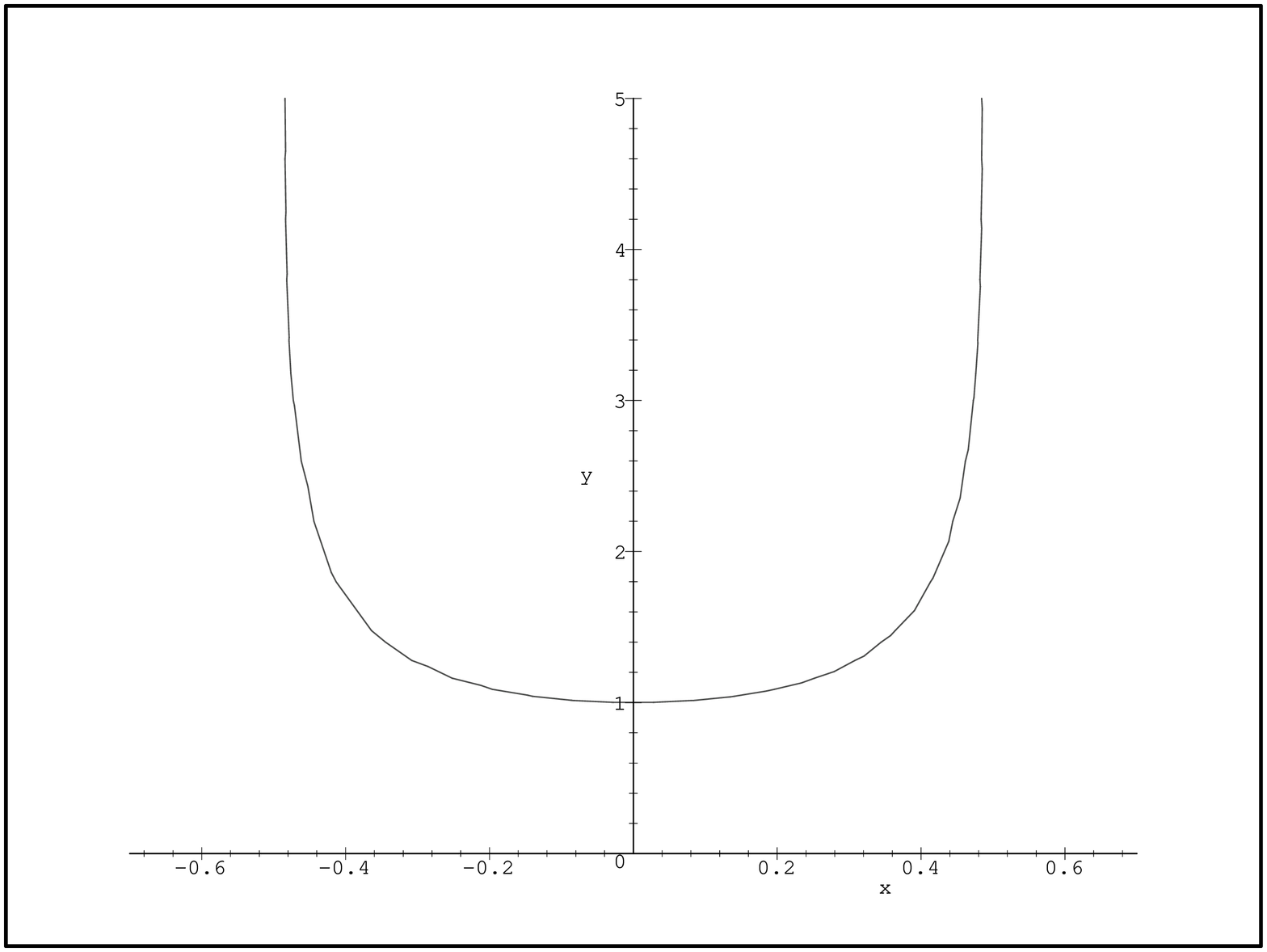}

Fig. 2.  Part of the A-polynomial curve corresponds to the curve $C$ in the image $i^*(\chi_{\uSL}(X)).$   
\end{center}

 We consider two surgery manifolds $X_1$ and $X_2.$ It follows from the computation in \cite{vu} that on each surgery manifold there is  a unique character and both characters have hyperbolic boundary holonomies. The  intersection between the curve $C$ and the curve $xy=1$ is the point $Q_1=(\alpha, 1/\alpha), \alpha\approx 0.4474073272,$ which corresponds to the representation $\rho_1$ on $X_1.$ The  intersection between the curve $C$ and the curve $x^2y=1$ is the point $Q_2=(\beta, 1/\beta^2), \beta\approx 0.4845486882,$ which corresponds to the representation $\rho_2$ on $X_2.$ 

It follows from \cite{vu} that $X_1$ is the homology sphere $\Sigma(2,3,11)$ and $GV(\rho_1|_X)=-\frac{50\pi^2}{33}$ 
and that  $X_2$ is the Seifert fibered space over $\S^2$ with three exceptional fibers of indices 2, 4 and 7 and  $GV(\rho_2|_X)=-\frac{9\pi^2}{7}$ . Now, using Theorem 2.1, we get the following identity:
$$GV(\rho_2|_X) - GV(\rho_1|_X)=4\int_\alpha^\beta \log y\frac{dx}{x}-\log x\frac{dy}{y}=-\frac{9\pi^2}{7}+\frac{50\pi^2}{33}=\frac{53\pi^2}{231}.$$
So we get the identity: $$\int_\alpha^\beta \log x\frac{dy}{y}-\log y\frac{dx}{x}=-\frac{53\pi^2}{924}.$$ Here the integral is taken over the curve $C$ from $Q_1$ to $Q_2.$
 
{\bf Remarks and questions:}
1) If the A-polynomial curve is of genus $0,$ then we can integrate the $1$-form $\log x\frac{dy}{y}-\log y\frac{dx}{x}$ in terms of the dilogarithm function. Unfortunately,  all the  knots that we consider here have the A-polynomial curves of genus  bigger than $1$ and we do not know how to write this kind of integrals in terms of the dilogarithm function.
  However we expect that one may do so since the volume of ideal simplexes in the $\uSL$ geometry  can be expressed in terms of the Roger dilogarithm function \cite{dupont}. 
  This would gives a new method for producing dilogarithm identities.

2) In \cite{boydvil1,boydvil2}, by writing the hyperbolic volume as sum of the Bloch-Wigner dilogarithms, the authors can express the Mahler's measure in terms of  Bloch-Wigner dilogarithm.
Our result here is in the same direction as \cite{boydvil1,boydvil2},  but the paths over which we integrate are in the real part of the A-polynomial curve whereas the integrals appearing in the Mahler's measure are taken over the imaginary part. 

3) Using the method in \cite{vu} one can find the exact value of the  Godbillon-Vey invariant of a representation on any Seifert fibered manifold as a rational multiple of $\pi^2.$
So whenever surgery on a knot gives a   Seifert fibered manifold, one would expect to find integration identity as we did above.

4) Can the value of the integrals above be predicted by results  from number theory?

\end{document}